\documentclass{amsart}
\usepackage{amssymb}
\usepackage{graphics}
\newtheorem{thm}{Theorem}[section]
\newtheorem{cor}[thm]{Corollary}
\newtheorem{lem}[thm]{Lemma}
\newtheorem{prop}[thm]{Proposition}
\theoremstyle{definition}
\newtheorem{defn}[thm]{Definition}
\theoremstyle{remark}
\newtheorem{rem}[thm]{Remark}
\numberwithin{equation}{section}
\newcommand{\abs}[1]{\left\vert#1\right\vert}

\newcommand{\Real}{\mathbb R}
\newcommand{\eps}{\varepsilon}
\newcommand{\To}{\longrightarrow}

\newcommand{\im}{\mathop{\textrm{Im}}\nolimits}
\newcommand{\imm}{\mathop{\textrm{Imm}}\nolimits}
\newcommand{\lgen}{\mathop{\textrm{Lgen}}\nolimits}
\newcommand{\proj}{\mathop{\textrm{Proj}}\nolimits}
\newcommand{\mono}{\mathop{\textrm{MONO}}\nolimits}
\newcommand{\imekv}[2]{\im(#1) \sim \im(#2)}
\begin{document}

\title{Regular homotopy classes of locally generic mappings}%
\author{Andr\'as Juh\'asz}%
\address{Department of Analysis, E\"otv\"os Lor\'and University,
P\'azm\'any P\'eter s\'et\'any 1/C, Budapest, Hungary 1117}%
\email{juhasz.6@dpg.hu}%

\thanks{Research partially supported by OTKA grant no. T037735}%
\subjclass{57R45; 58K30; 57R42}%
\keywords{immersion, locally generic mapping, regular homotopy, cross-cap singularity}%

%\date{\today}%
%\dedicatory{}%
%\commby{}%
% ----------------------------------------------------------------
\begin{abstract}
In this paper we generalize the notion of regular homotopy of
immersions of a closed connected $n$-manifold into $\Real^{2n-1}$
to locally generic mappings. The main result is that if $n=2$ then
two mappings with singularities are regularly homotopic if and
only if they have the same number of cross-cap (or
Whitney-umbrella) singularities. As an application, we get a
description of the path-components of the space of those
immersions of a surface into $\Real^4$ whose projections into
$\Real^3$ are locally generic.
\end{abstract}
\maketitle
% ----------------------------------------------------------------
\section{Introduction}

Our work was motivated by the paper of U. Pinkall \cite{Pinkall},
which classifies immersions of compact surfaces into $\Real^3$ up
to regular homotopy, allowing diffeomorphisms of the source
manifold $M^2$. So two immersions $ f,g \colon M^2 \To R^3 $ are
considered equivalent if there is a diffeomorphism $\varphi$ of
$M^2$ such that $f = g \circ \varphi$. This notion is different
from regular homotopy, yielding an interesting classification of
immersed surfaces using the Arf invariant. That paper also gives
generators for the abelian semigroup of immersed surfaces with the
connected sum operation. Professor Andr\'as Sz\H{u}cs asked me
what happens with Pinkall's classification if we allow cross-cap
(also called Whitney-umbrella) singularities. The notion of
regular homotopy has to be revised and, unlike for immersions, for
singular maps it turns out that all natural definitions are
equivalent. In fact, we prove that for singular mappings (i.e. not
immersions) of a closed connected surface the number of cross-caps
totally determines the regular homotopy class, diffeomorphisms of
the source manifold are not needed. Thus the approach of U.
Pinkall and the classical regular homotopy classification give the
same result for singular maps.

In the final section of our paper we present an application of the
above result to the study of the path components of the space of
those immersions of a closed connected surface into $\Real^4$
whose projections into $\Real^3$ are locally generic, i.e. may
have cross-cap singularities.

It is a natural question if Theorem \ref{thm:1} generalizes to
higher dimensions, for locally generic maps of a closed n-manifold
$M^n$ into $\Real^{2n-1}$. Our methods of proof for Theorem
\ref{thm:1} do not seem to work if $n > 2$ since they rely heavily
on the results of surface topology. However, I could prove the
general result in the case when $n >3 $ and $M^n$ is 2-connected.
I will publish this in a separate paper.  I want to emphasize that
the results of Section 3 easily generalize for any closed manifold
$M^n$ provided that the generalization of Theorem \ref{thm:1}
holds true for $M^n$.

I would like to take this opportunity to express my gratitude to
Professor Andr\'as Sz\H{u}cs for drawing my attention to this
problem and for his constant support and encouragement. I would
also like to thank Professor Bal\'azs Csik\'os who read the first
version of this paper and suggested several improvements.

\section{The main result}

Let us start by introducing the notion of a locally generic
mapping of a closed $n$-manifold $M^n$ into a $(2n-1)$-manifold
$N^{2n-1}$ for $n \ge 2$.

\begin{defn}
$f \colon M^n \To N^{2n-1}$ is called \emph{locally generic} if it
is an immersion except for cross-cap singularities. The set of
singular points of $f$ in $M^n$ is denoted by $S(f)$.
($\abs{S(f)}<\infty$ because $M^n$ is compact.) Let us denote by
$\lgen(M^n, N^{2n-1})$ the subspace of locally generic mappings in
$C^{\infty}(M^n, N^{2n-1})$ endowed with the $C^{\infty}$
topology.
\end{defn}

Recall that a map $f \colon M^n \To N^{2n-1}$ is called generic
(or stable), if it is an immersion with normal crossings except in
a finite set of points, moreover the singular points of $f$ are
non-multiple cross-cap points. This explains our terminology.
Whitney \cite{Whitney} proved that the set of stable maps is dense
open in $C^{\infty}(M^n, \Real ^{2n-1})$ with respect to the
$C^{\infty}$-topology. Studying the double point set of a generic
mapping close to the locally generic mapping $f$ in the
$C^{\infty}$-topology, one can easily verify that for $M^n$ closed
$\abs{S(f)}$ is an even integer. (For the cross-caps are precisely
endpoints of double-point curves.)

\begin{defn}
\label{defn:1} Two locally generic mappings $f, g \colon M^n \To
N^{2n-1}$ are called \emph{regularly homotopic} ( denoted by $f
\sim g$ ) if there is a smooth mapping $H \colon M^n \times [0,1]
\To N^{2n-1}$ such that $H_t$ is locally generic for each $t \in
[0,1]$, moreover $H_0 = f$ and $H_1 = g$. Here $H_t(x) = H(x,t)$
for $x \in M^n$ and $t \in [0,1]$.
\end{defn}

The following definition will be especially useful in the case
$n=2$.

\begin{defn}
Two locally generic mappings $f, g \colon M^n \To N^{2n-1}$ are
called \emph{image-homotopic} if there is a diffeomorphism
$\varphi$ of $M^n$ such that $f \circ \varphi$ is regularly
homotopic to $g$. We denote this by $\imekv{f}{g}$.
\end{defn}

\begin{prop} \label{prop:1}
If $f \sim g$ or just $\imekv{f}{g}$, then $\abs{S(f)} =
\abs{S(g)}$. In fact, if $H \colon M^n \times [0,1] \To N^{2n-1}$
is a regular homotopy, then $\abs{S(H_t)}=\abs{S(H_0)}$ for every
$t \in [0,1] $.
\end{prop}

\begin{proof}
If $\varphi$ is a diffeomorphism of $M^n$, then $S(f) = \varphi(
S(f \circ \varphi))$, so $\abs{S(f)} = \abs{S(f \circ \varphi)}$.
Thus we can suppose that $f \sim g$. From the definition of
stability it is clear that every locally generic mapping $h$ has a
neighborhood $U_h$ in the Whitney $C^{\infty}$ topology such that
for every $h' \in U_h$ we have $\abs{S(h')} = \abs{S(h)}$ (because
every $p \in M^n$ has a neighborhood $U$ such that $h|U$ is
equivalent to $h'|U$ and since $M^n$ is closed, see \cite[p.
72]{sing}). Thus the function $\abs{S(.)} \colon \lgen(M^n,
N^{2n-1}) \To \mathbb Z$ is locally constant. So if $H$ is a
regular homotopy connecting $f$ and $g$, this implies that $H_t$
is a continuous path in $\lgen(M^n, N^{2n-1})$, along which
$\abs{S(H_t)}$ is constant.
\end{proof}

\begin{rem}
We have defined the notions of regular homotopy and image-homotopy
between two locally generic maps $f$ and $g$. Proposition
\ref{prop:1} implies that if $f$ and $g$ are immersions then they
are regularly homotopic, resp. image-homotopic as immersions iff
they are those as locally generic maps.
\end{rem}

If $H \colon M^2 \times [0,1] \To \Real^3$ is a regular homotopy
between two locally generic mappings $f = H_0$ and $g = H_1$ and
$k = \abs{S(f)} = \abs{S(g)}$, then there exists curves $\gamma_1,
\dots \gamma_k \colon [0,1] \To M^2$ such that $S(H_t) = \left\{
\, \gamma_j(t) \colon 1 \le j \le k \, \right\}$ for every $t \in
[0,1]$. We define the bijection $i_H \colon S(f) \To S(g)$ the
following way : for $1 \le j \le k$ let $i_H(\gamma_j(0)) =
\gamma_j(1)$.

Now we can state the main result of this paper yielding a converse
of Proposition \ref{prop:1} for singular mappings.

\begin{thm}
\label{thm:1} Let $M^2$ be a closed connected surface and suppose
that $f,g \colon M^2 \To \Real^3$ are locally generic mappings
with $\abs{S(f)} = \abs{S(g)} > 0$. Then $f \sim g$. Moreover for
any bijection $i \colon S(f) \To S(g)$ there exists a regular
homotopy $H$ connecting $f$ and $g$ such that $i_H = i$.
\end{thm}

Let us now list a few interesting corollaries of this. A
surprising consequence of Theorem \ref{thm:1} is the following: If
$U$ is the standard locally generic mapping of $S^2$ into
$\Real^3$ with two cross-cap points then for any two immersions
$f, g \colon M^2 \To \Real^3$ the connected sums $f\#U$ and $g\#U$
become regularly homotopic as locally generic maps !

Another consequence of our theorem is that, given a closed
connected surface $M^2$, we can easily produce a full list of
representatives of all regular homotopy classes of locally generic
maps $M^2 \To \Real^3$. (We only do this for singular maps, for
immersions see \cite{Pinkall}.) First suppose that $M^2$ is
orientable. If we denote by $i_M$ the standard embedding of $M^2$
into $\Real^3$ then $$i_M  \# \underbrace{U \# \dots \# U}_n$$ is
a representative for the class of locally generic maps with $2n
>0$ singular points. Now suppose that $M^2$ is a non-orientable
surface of genus $g$. We denote the Boy surface by $B$. Now
$$\underbrace{B \# \dots \# B}_g \# \underbrace{U \# \dots \#
U}_n$$ is a locally generic mapping of $M^2$ into $\mathbb{R}^3$
with $2n$ cross-cap points.

\begin{figure}
\includegraphics{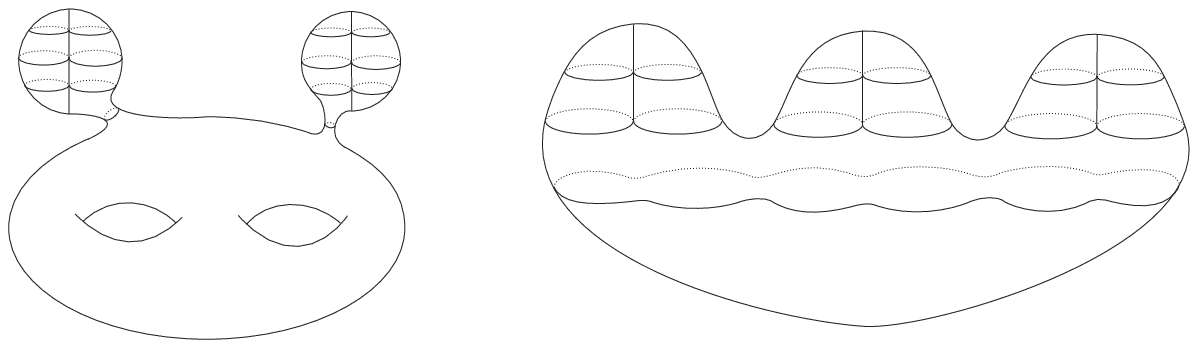}
\caption{} \label{fig1}
\end{figure}

Perhaps the following construction can be visualized more easily:
Let us denote by $V$ the well-known locally generic mapping of
$\Real P^2$ into $\Real^3$ having two singular points ($V \sim B
\# U$). Then any singular locally generic mapping is regularly
homotopic to one of the form $i_N \# U \# \dots \# U \# V \# \dots
\# V$, where $N^2$ is orientable. (The left side of Figure
\ref{fig1} depicts $A_2 \# U \# U$, where $A_2$ is the orientable
surface of genus $2$. The right side of Figure \ref{fig1}
illustrates $B\# B \# B \# U \# U \# U \sim V \# V \# V$.)

Pinkall determined the abelian semigroup $H$ of immersed surfaces
in $\Real^3$ with the connected sum operation (see
\cite{Pinkall}). If we consider the extended semigroup
$\widetilde{H}$ of locally generic surfaces, then only $U$ is
needed as a new generator with the following new relations (using
Pinkall's notation): $S \# U  = T \# U$ and $B \# U = \overline{B}
\# U$ .

Note that $\widetilde{H} \setminus H$ is a sub-semigroup of
$\widetilde{H}$ and it easily follows from Theorem \ref{thm:1}
that it is isomorphic to $J \oplus \mathbb{Z}_+$, where $J$
denotes the semigroup of (closed connected) surfaces. As a
corollary we may conclude that the Grothendieck group of
$\widetilde{H}$ is isomorphic to $\mathbb{Z} \oplus \mathbb{Z}$.

\section{Proof of the main result}

The purpose of this section is to prove Theorem \ref{thm:1}.\
First we recall the classification of immersions of an arbitrary
two-dimensional manifold $F^2$ into $\Real^3$ using Hirsch-Smale
theory.

\begin{thm}
\label{thm:3} There is a 1-1 correspondence between the regular
homotopy classes of immersions of a surface $F^2$ into $\Real^3$
denoted by $\imm(F^2, \Real^3)$ and $H^1(F^2 ; \mathbb{Z}_2)$.
\end{thm}

\begin{proof}
By Hirsch (see \cite{Hirsch}) there is a  weak homotopy
equivalence between the space $\imm(F^2, \Real^3)$ and
$\Gamma(\mu)$. Here $\Gamma(\mu)$ denotes the space of sections of
the vector bundle $\mu = \mono(TF^2,R^3)$ over $F^2$ whose fiber
over $p \in F^2$ consists of all linear injections from $T_pF$ to
$\Real^3$. Thus there is a bijection between the regular homotopy
classes $\pi_0(\imm(F^2, \Real^3))$ and $\pi_0(\Gamma(\mu))$. Fix
an arbitrary Riemannian metric on $F^2$ and let $\mu'$ be the
bundle over $F^2$ whose general fiber over $p$ is the space of
orthogonal injections of $T_pF$ into $\Real^3$. Then the inclusion
of $\mu'$ into $\mu$ is a fiber homotopy equivalence (see \cite[p.
426]{Pinkall}), thus $\pi_0(\Gamma(\mu))=\pi_0(\Gamma(\mu'))$.
Fixing a section $s \in \Gamma(\mu')$ every section $t \in
\Gamma(\mu')$ can be obtained by the action of a unique element of
$C(F^2, SO(3))$ on $s$. Thus $\Gamma(\mu')$ is homeomorphic to
$C(F^2, SO(3))$, yielding $\pi_0(\Gamma(\mu))=[F^2, SO(3)]$. Since
$SO(3)$ is homeomorphic to $RP^3$, it follows from obstruction
theory that
\[[F^2, SO(3)]=[F^2, RP^3]=[F^2, RP^{\infty}]=[F^2,
K(\mathbb{Z}_2,1)]=H^1(F^2; \mathbb{Z}_2),\] where
$K(\mathbb{Z}_2,1) = RP^{\infty}$ is an Eilenberg-MacLane space.
\end{proof}

From now on $M^2$ denotes the closed connected surface mentioned
in the statement of Theorem \ref{thm:1}. If $f$ and $g$ are
locally generic mappings of $M^2$ into $\Real^3$ with
$\abs{S(f)}=\abs{S(g)}$, then according to the lemma of
homogeneity there exists a diffeotopy $\{\varphi_t \colon t \in
[0,1] \}$ of $M^2$ such that $\varphi_0 = id_{M^2}$ and
$\varphi_1(S(g))=S(f)$. Since $f \circ \varphi_t$ provides a
regular homotopy between $f$ and $f \circ \varphi_1$, it is
sufficient to prove Theorem \ref{thm:1} in the case $S(f)=S(g)$.
Let $S(f) = S(g)= \{\, p_1, \dots , p_k \,\}$, where $k = \abs
{S(f)} > 0$ is an even integer. For each $p_i$ choose a
sufficiently small open neighborhood $D_i$ diffeomorphic to an
open 2-disc such that $f|D_i$ has the canonical form $(x_1^2, x_2,
x_1x_2)$ in an appropriate pair of local coordinate-systems
centered at $p_i$ and $f(p_i)$. Similarly $g|D_i$ should have the
same canonical form in another pair of local coordinate-systems.
Assume moreover that the discs $D_1, \dots, D_k$ are pairwise
disjoint. Denote by $A$ the disjoint union $\coprod_{i=1}^k D_i$,
then $F^2 = M^2 \setminus A$ is a two-manifold with boundary.

\begin{lem} \label{lem:1}
Suppose that $f$ and $g$ are locally generic mappings of the
surface $M^2$ into $\Real^3$ such that $S(f)=S(g)$. Choose open
discs $D_1, \dots, D_k \subset M^2$ centered at the points of
$S(f)$ as above. Define $F^2 = M^2 \setminus \coprod_{i=1}^k D_i$.
Then there exists a diffeomorphism $d$ of the pair $(M^2,F^2)$
such that the immersions $(f|F^2) \circ d$ and $g|F^2$ are
regularly homotopic and $d$ permutes $S(f)$. Moreover there is a
diffeotopy $d_t$ of $M^2$ with $d_0 = id_{M^2}$ and $d_1 = d$.
\end{lem}

\begin{proof}
According to Theorem \ref{thm:3} the regular homotopy classes of
$f|F^2$ and $g|F^2$ correspond to cohomology classes $\alpha,
\beta \in H^1(F^2 ; \mathbb{Z}_2)$. (These will be shown to be
non-zero later.) We construct a diffeomorphism $d$ of the pair
$(M^2, F^2)$ such that for the induced automorphism $d^*$ of
$H^1(F^2 ; \mathbb{Z}_2)$ it holds that $d^*(\alpha) = \beta$.
Using Theorem \ref{thm:3} again this gives the required result
$(f|F^2) \circ d \sim g|F^2$.

We first note that \begin{equation} \label{eqn:0} H_1(F^2 ;
\mathbb{Z}_2) = H_1(M^2; \mathbb{Z}_2) \oplus
\overbrace{\mathbb{Z}_2 \oplus \dots \oplus \mathbb{Z}_2}^{k-1}\,,
\end{equation} as can be seen from the exact sequence of the pair $(M^2, F^2)$.
(Recall that $k = \abs{S(f)}$.) For each $i$, $1 \le i \le k$
choose an embedded curve $\gamma_i$ in $F^2$ around $D_i$. Denote
its homology class by $[\gamma_i] = c_i$. The classes $c_1, \dots,
c_{k-1}$ can be chosen for the generators of the $(k-1)$
$\mathbb{Z}_2$ summands in \ref{eqn:0}. (Note that $c_1 + \dots +
c_k =0$ since $\partial F^2 = \gamma_1 \cup \dots \cup \gamma_k.$)
According to Pinkall \cite{Pinkall} we have that
$\langle\alpha,c_i\rangle= 1$ and $\langle\beta,c_i\rangle=1$ for
every $i \in \{\,1, \dots, k \,\}$, since $\gamma_i$ has a
neighborhood homeomorphic to $S^1 \times [0,1]$ which is mapped by
$f$ and also by $g$ into a "figure eight$\times [0,1]$". Now we
have two cases according to the orientability of $M^2$.

If $M^2$ is orientable of genus $g$ we denote the standard
generators of $H_1(M^2 ; \mathbb{Z}_2)$ by $a_1, b_1, \dots, a_g,
b_g$, and choose embedded curves $\varphi_1, \psi_1 \dots
\varphi_g, \psi_g$ in $F^2$ representing them. Define
\[H_a = \{\,1 \le i \le g \colon \langle \alpha , a_i \rangle \neq
\langle \beta, a_i \rangle \,\},\] and \[H_b = \{\,1 \le i \le g
\colon \langle \alpha , b_i \rangle \neq \langle \beta, b_i
\rangle \,\}.\] There is a simple (i.e. embedded) closed curve
$\delta$ in $F^2$ that for each $i$, $1 \le i \le g$ intersects
transversally in one point the curve $\varphi_i$ if $i \in H_a$
and is disjoint from $\varphi_i$ if $i \not\in H_a$, moreover
$\delta$ intersects transversally in one point $\psi_i$ if $i \in
H_b$ and is disjoint from $\psi_i$ if $i \not\in H_b$. Note that
the homology class of $\delta$ will be
\begin{equation} \label{eqn:1} \left(\sum_{i \in H_a} b_i \right) + \left(\sum_{i
\in H_b} a_i \right).\end{equation} Such a $\delta$ exists because
in $H_1(M^2; \mathbb{Z}_2)$ any class can be represented by a
simple curve, and a simple curve representing the class
\ref{eqn:1} can be arranged to be transversal to all the curves
$\varphi_j$ and $\psi_j$ for $j = 1, \dots, k$ and to intersect
each of them at most in one point. Now choose two points on
$\delta$ very close to each other so that none of the curves
$\varphi_i$ and $\psi_i$ for $1 \le i \le g$ intersects the
shorter arc $\delta^*$ between them.  Thus the following
equalities hold :
\begin{equation} \label{eqn:2}
\abs{\delta \cap \varphi_i} = \left\{
\begin{array}{ll}
1 & \mbox{if $i \in H_a$}\\
0 & \mbox{if $i \not\in H_a$}
\end{array}
\right.
\end{equation}
\begin{displaymath}
\abs{\delta \cap \psi_i} = \left\{
\begin{array}{ll}
1 & \mbox{for $i \in H_b$}\\
0 & \mbox{for $i \not\in H_b$}
\end{array}
\right.
\end{displaymath}
\begin{eqnarray}
\delta^* \cap \varphi_i = \delta^* \cap \psi_i = \emptyset &
\mbox{for $1 \le i \le g$.} \nonumber
\end{eqnarray}
Modify the arc $\delta^*$ ( by a homology ) so that it goes
through the center of $D_1$ and $D_2$ avoiding the curves
$\varphi_i$ and $\psi_i$ for $i = 1, \dots, g$ as well as the
discs $D_i$ for $i = 3, \dots, k$. This can be done since $M^2
\setminus \bigcup_{i=1}^{2g} (\varphi_i \cup \psi_i)$ is
path-connected. From now on we will denote this new simple curve
on $M^2$ by $\delta$ (see Figure \ref{fig2}). Note that the
equalities \ref{eqn:2} still hold. Now choose a tubular
neighborhood $T$ of $\delta$ such that $D_1, D_2 \subset T$ and
for $i \in H_a$ the curve $\varphi_i$ and for $j \in H_b$ the
curve $\psi_j$ intersects $T$ in a line segment. (See the left
side of Figure \ref{fig3}.) We also select a slightly wider
tubular neighborhood $T' \supset T$. Define $d$ on $T$ to be a
rotation of $T = S^1 \times [-1,1]$ by $180^{\circ}$ interchanging
$D_1$ and $D_2$ and also interchanging $p_1$ and $p_2$. The
diffeomorphism $d$ acts identically on $M^2 \setminus T'$. On $T'
\setminus T$, which is homeomorphic to $S^1 \times
([-2,2]\setminus [-1,1])$, define $d$ as the rotation of $S^1
\times \{s\}$ by $(2- \abs{s}) \times 180^{\circ}$ for $s \in
[-2,2] \setminus [-1,1]$ (see the right side of Figure
\ref{fig3}). The diffeomorphism $d$ is diffeotopic to $id_{M^2}$ :
construct $d_t \colon M^2 \To M^2$ similarly to $d$, just take a
rotation by $180^{\circ} \cdot t$ instead of rotating by
$180^{\circ}$ such that $d_t(p_1) \in \delta^*$ for every $t \in
[0,1]$.
\begin{figure}
\includegraphics{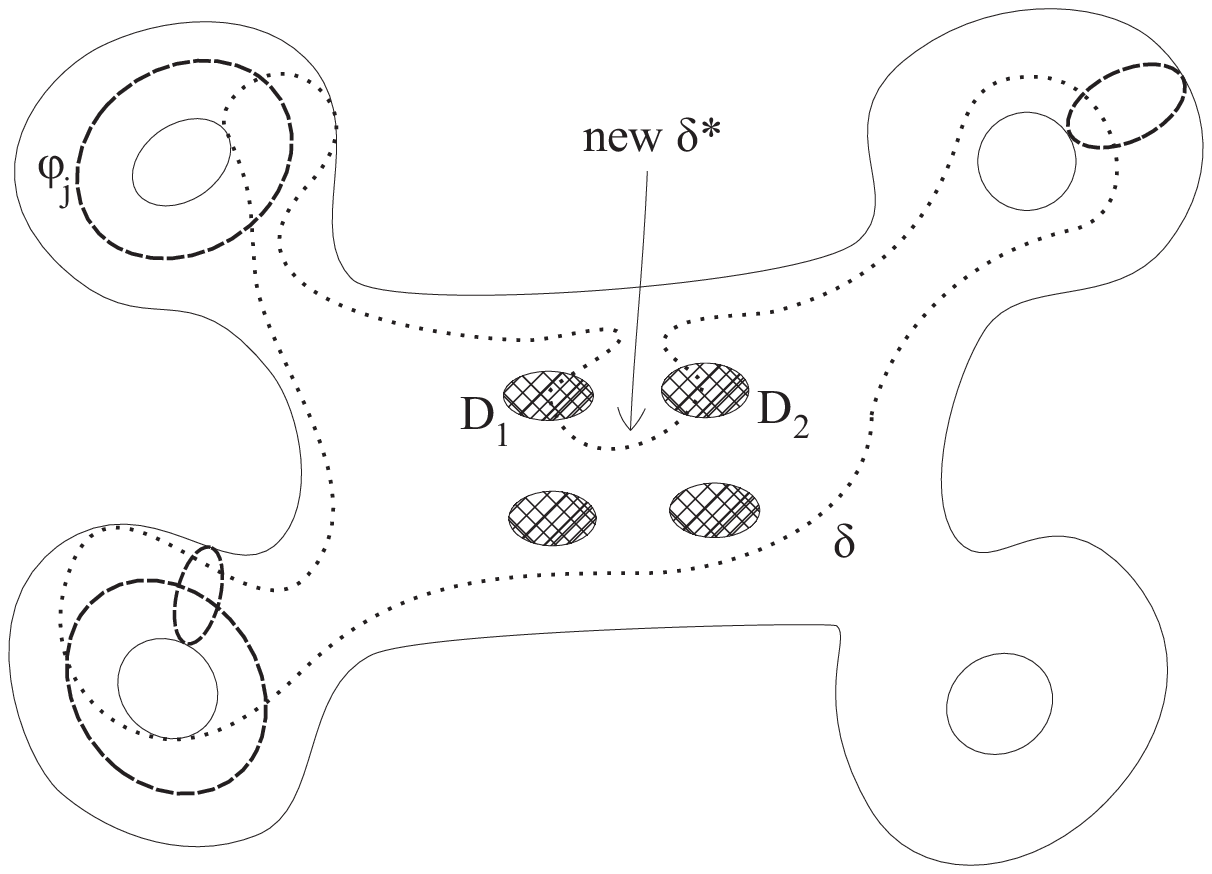}
\caption{} \label{fig2}
\end{figure}

\begin{figure}
\includegraphics{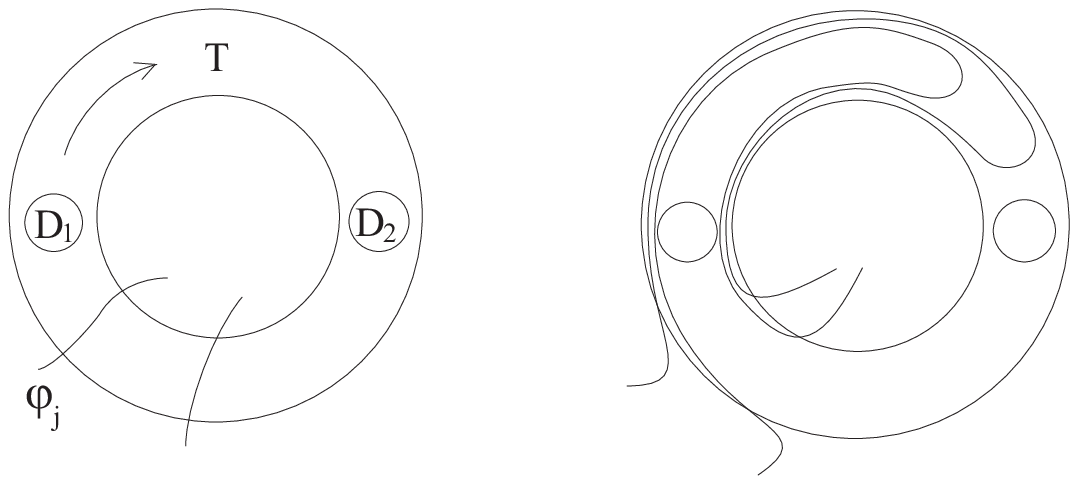}
\caption{} \label{fig3}
\end{figure}

For $i \not\in H_a$ it is clear that $d$ is the identity on the
image of $\varphi_i$, thus $d_*(a_i)=a_i$. On the other hand if $i
\in H_a$, then $d \circ \varphi_i$ is homologous to the connected
sum of $\varphi_i$ and $\gamma_1$ (surrounding $D_1$), thus
$d_*(a_i) = a_i + c_1$. Similarly for $i \not\in H_b$ we have
$d_*(b_i)=b_i$, and if $i \in H_b$ then $d_*(b_i) = b_i + c_1$.
Finally it holds that $d_*(c_2)=c_1$, $d_*(c_1)=c_2$ and if $i>2$
then $d_*(c_i)=c_i$. (If $k=2$ then $d_*(c_1) = c_2 = c_1$.) This
can be verified by looking at the action of $d$ on the curves
$\gamma_1, \dots, \gamma_k$. Since $d_*$ permutes the generators
$c_1, \dots c_{k-1}$, for $1 \le i \le k-1$ we have $\langle
\alpha, d_*(c_i) \rangle = 1$. (Recall that $\langle \alpha, c_i
\rangle = 1$ and $\langle \beta, c_i \rangle = 1$ for every $i$.)
Thus $\langle \alpha, d_*(c_i) \rangle = \langle \beta, c_i
\rangle$. By the choice of $H_a$ we see that for $i \in H_a$ it
holds that $\langle \alpha, d_*(a_i) \rangle = \langle \alpha, a_i
+ c_1 \rangle = \langle \alpha, a_i \rangle + 1 = \langle \beta,
a_i \rangle$ and for $i \not\in H_a$ we have that $\langle \alpha,
d_*(a_i) \rangle = \langle \alpha, a_i \rangle = \langle \beta,
a_i \rangle$. A similar argument holds for $b_1, \dots, b_g$.
Since $a_1, b_1, \dots, a_g, b_g$ and $c_1, \dots, c_{k-1}$ form a
basis of $H_1(F^2; \mathbb{Z}_2)$, we have shown that $\langle
d^*(\alpha), x \rangle = \langle \alpha, d_*(x) \rangle = \langle
\beta, x \rangle$ for every $x \in H_1(F^2; \mathbb{Z}_2)$. Thus
$d^*(\alpha)=\beta$. Hence $(f \circ d )| F^2$ is regularly
homotopic to $g|F^2$ as required. Also $d$ satisfies
$d(S(f))=S(f)$.

Now suppose that $M^2$ is a non-orientable surface of genus $g$,
i.e. a sphere with $g$ Moebius bands. Choose a curve $\varphi_i$
in $F^2$ on the $i$-th Moebius band representing its homology
generator for $1 \le i \le g$. Then $b_1 = [\varphi_1] , \dots,
b_g = [\varphi_g]$ together with $c_1, \dots, c_{k-1}$ is the
standard basis of $H_1(F^2; \mathbb{Z}_2)$. Analogously to the
orientable case let
\[H = \{\,1 \le i \le g \colon \langle \alpha , b_i \rangle \neq
\langle \beta, b_i \rangle \,\},\] and similarly it is sufficient
to construct a diffeomorphism $d$ of the pair $(M^2, F^2)$ with
$d_*(c_i)=c_i$ for $1 \le i \le k$, $d_*(b_i)=b_i + c_1$ for $i
\in H$ and $d_*(b_i) = b_i$ for $i \not\in H$. It is enough to
show that for any fix $1 \le j \le g$ we can find a diffeomorphism
$d_j$ of the pair $(M^2, F^2)$ such that $d_{j*}(b_j) = b_j + c_1$
and that $d_{j*}$ is identical on every other homology-generator.
(Then $\prod_{j \in H} d_j$ is a good choice for $d$.)

For this end modify a small arc of $\varphi_j$ using a homology
such that it still lies in $F^2$ but gets close to $D_1$ and
remains disjoint from all the other $\varphi_i$ for $i \neq j$.
(We shall call this modified curve $\varphi_j$ also.) This is
possible since $M^2 \setminus \bigcup_{i=1}^g \varphi_i$ is
path-connected. Denote by $T$ a tubular neighborhood of
$\varphi_j$ in $M^2$ containing $D_1$ and disjoint from $D_i$ if
$i > 1$ and from $\varphi_i$ if $i \neq j$. (See Figure
\ref{fig4}.) Then $T$ is homeomorphic to the Moebius band. Also
choose a slightly larger tubular neighborhood $T'$ of $T$ with
similar properties. Now think of $T$ as a rectangle with the
vertical sides identified in the opposite direction and with $D_1$
in its center. Let $d_j$ be the reflection of the rectangle $T$
into its horizontal central axis going through the center of $D_1$
which is $p_1$. Then $d_j$ induces an orientation-preserving
diffeomorphism ( a rotation ) of $\partial T = S^1$, which can be
extended to $T' \setminus T = S^1 \times [0,1]$ being identical on
$\partial T' = S^1 \times \{1\}$ as we have already seen. Finally
$d_j$ is identical on $M^2 \setminus T'$. This $d_j$ maps the
curve $\varphi_j$ (which is the horizontal central line in the
rectangle except that it avoids $D_1$ (see Figure \ref{fig4})) to
a curve homologous to the connected sum of $\varphi_j$ and
$\gamma_1$, thus $d_{j*}(b_j) = b_j + c_1$. Since $\varphi_i$  for
$i \neq j$ and $\gamma_i$ for $i
> 1$ are fixed by $d_j$, it satisfies the required conditions.
Concerning $c_1$ we have $d_{j*}(c_1) = - c_1 = c_1$ since we are
working with mod $2$ coefficients.
\begin{figure}
\includegraphics{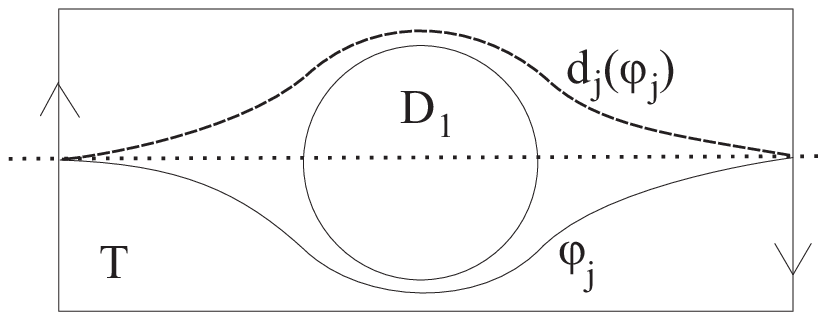}
\caption{} \label{fig4}
\end{figure}

The diffeomorphism $d_j$ of $M^2$ is diffeotopic to $id_{M^2}$:
Think of the Moebius band $T$ as the factor space $S^1 \times
[0,1] / _{p \times \{1\} = (-p) \times \{1\}}$ (thus we identify
the opposite points of one boundary component of an annulus). In
this model define $(d_j)_t$ on $T$ as the diffeomorphism induced
by the rotation of the annulus by $180^{\circ} \cdot t$ degrees.
On $T \setminus T'$ define $(d_j)_t$ as before. Finally on $M^2
\setminus T'$ the diffeomorphism $(d_j)_t$ is the identity
mapping. Then $d_j = (d_j)_1$ and $(d_j)_0 = id_{M^2}$ as
required.
\end{proof}

As a consequence of the above proof we obtain the following
proposition (for the definition of the mapping class group see
\cite{Stillwell}):

\begin{cor}
Suppose that $F^2$ is a surface of genus $g$ with $k>1$ boundary
components, where k is even. We denote, like as before, the
homology classes represented by the boundary components of $F^2$
in $H_1(F^2; \mathbb{Z}_2)$ by $ c_1, \dots, c_k$. The mapping
class group $M(F^2)$ of $F^2$ acts on the set $S = \{\,\alpha  \in
H^1(F^2, \mathbb{Z}_2) \colon \langle \alpha,c_i \rangle = 1, \, 1
\le i \le k-1 \,\}$. (If $\alpha \in S$ then $\langle \alpha, c_k
\rangle = \langle \alpha, c_1 + \dots + c_{k-1} \rangle = 1$ since
$k$ is even.) If $M^2$ is a closed surface of genus $g$ then there
is a homomorphism $m \colon M(F^2) \To M(M^2)$ obtained by
"filling in the holes". Then $\ker(m)$ acts transitively on $S$.
\end{cor}

\begin{lem}
\label{lem:2} Let $g$ and $h$ be locally generic mappings of $M^2$
into $\Real^3$ such that $S(g) = S(h)$ and $g|F^2 \sim h|F^2$,
where $F^2$ is the complement of a small open neighborhood of the
common singular set. Then $g \sim h$.
\end{lem}

\begin{proof}
Recall that $A = \coprod_{i=1}^k D_i$. Since for each $i$ it holds
that $h|D_i$ and $g|D_i$ have canonical forms in appropriate
coordinate-systems, there is a regular homotopy $H$ between $h$
and a locally generic mapping $\widetilde{h}$ such that for each
$t \in [0,1]$ we have $S(H_t)=S(h)$ and that $\widetilde{h}|A =
g|A$. Thus $H|(F^2 \times [0,1])$ is a regular homotopy between
the immersions $h|F^2$ and $\widetilde{h}|F^2$ showing that
$\widetilde{h}|F^2 \sim g|F^2$. So we can suppose that $g|A =
h|A$.

Let $H$ be a regular homotopy between $g|F^2$ and $h|F^2$. For
every $i = 1, \dots, k$ fix a smaller concentric closed disc $B_i$
in $D_i$ (hence $p_i \in B_i \subset D_i$). Finally set $F_l = F^2
\cup \left(\bigcup_{i=1}^{l} D_i \right)$. We will define
recursively a sequence of regular homotopies $H^l \colon F_l
\times [0,1] \To \Real^3$ connecting $g_l = g|F_l$ and $h_l =
h|F_l$ for $l=0,\dots,k$ with the property $H^0 = H$. Suppose that
we have constructed $H^l$ for $l < j$. Let $q$ be a point in
$\partial D_j$. For $t \in [0,1]$ there is a one-parameter family
of elements $M_t \in GL(\Real, 3)$ with $M_t \circ d_qH^{j-1}_t =
d_qH^{j-1}_0$ and a vector $v_t \in \Real^3$ with
$M_t(H^{j-1}_t(q)) + v_t = H^{j-1}_0(q)$. Here $d_qH^{j-1}_t$
denotes the differential of the mapping $H^{j-1}_t$ at the point
$q$. Now define $H^j_t|F_{j-1}$ to be equal to $M_t \circ
H^{j-1}_t + v_t$. Since $g(q) = h(q)$ and $d_qg =d_qh$, we have
that $M_1 = id_{\Real^3}$ and $v_1 = 0$, thus $H^j_1|F_{j-1} =
h|F_{j-1}$. With this transformation of $H^{j-1}$ we have achieved
that $H^j_t(q) = H^j_0(q)$ and $d_qH^j_t = d_qH^j_0$ for every $t
\in [0,1]$. Let $U_q$ be a small closed neighborhood of $q$ in
$F_{j-1}$ diffeomorphic to a closed 2-disc ($q \in \partial U_q$).
Using a standard argument of S. Smale we can suppose that the
homotopy $H^j$ is kept fixed on a neighborhood of $U_q$ (see
Hirsch \cite{Hirsch}, Lemma 2.5 on page 249). On $B_j$ define for
every $t \in [0,1]$ the mapping $H^j_t|B_j = H^j_0|B_j = g|B_j$.

\begin{figure}
\includegraphics{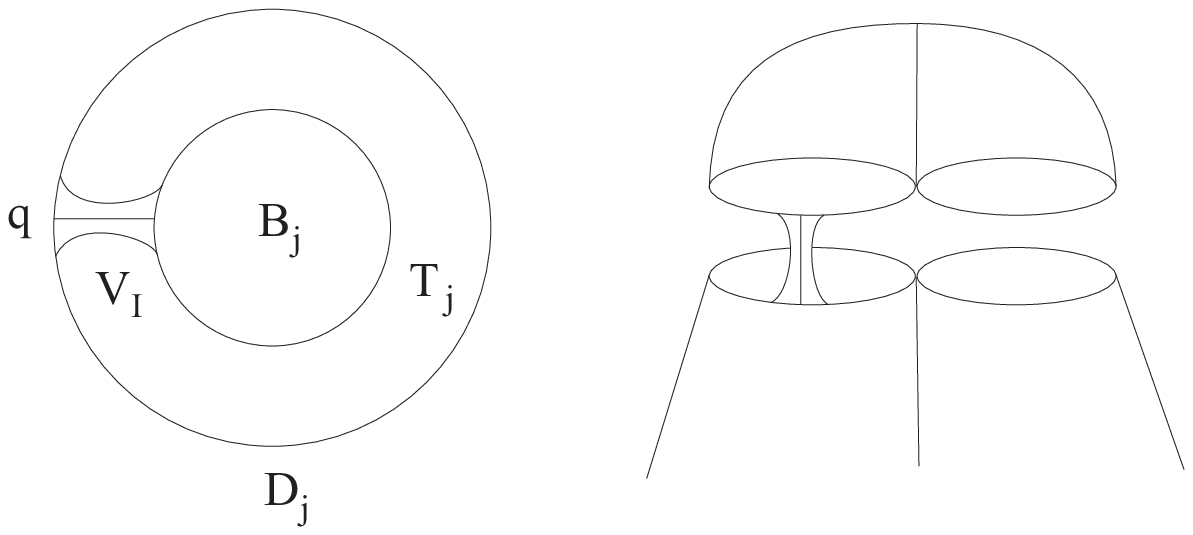}
\caption{} \label{fig5}
\end{figure}

Denote the closure of the annulus $D_j \setminus B_j$ by $T_j$,
and let $I$ be a radial line in $T_j$ containing $q$ (see Figure
\ref{fig5}). A tubular neighborhood $V_I \subset T_j$ of $I$ is
obtained by taking $V_I = (U_q \cap
\partial D_j) \times I$. Let $H^j_t|V_I = g|V_I$ for every $t \in
[0,1]$. This is possible since $H^j$ is fix on $U_q$. We only have
to define $H^j$ on the closed two-cell $C_i = T_i \setminus
\text{Int}( V_I)$. For this purpose we will use Smale's lemma (see
Theorem 1.1 on page 245 of \cite{Hirsch} or Theorem 2.1 in
\cite{Smale}), which intuitively states the following: If we are
given an immersed disk $D^k$ in $\Real^n$ such that $k < n$ and we
deform the boundary of the disk and the normal derivatives along
the boundary, then we can deform the whole disk at the same time
so as to induce the given deformation on the boundary and normal
derivatives. Since $H^j$ is already defined on $\partial C_i$
along with derivatives normal to $\partial C_i$ we can use Smale's
lemma to get a regular homotopy $G$ on $C_j$ with $G_0 = g|C_j$
and $G$ prescribed along the boundary. Finally $G_1 \sim h|C_j$
since they coincide on $\partial C_j$ and the obstruction is an
element of $\pi_2(SO(3)) = 0$. Putting this homotopy after $G$ we
obtain the desired homotopy $H^j|C_j$. Thus we have constructed
$H^j$ on the whole manifold $F_j$. This shows that $H^k$ is a
regular homotopy connecting $g$ and $h$.
% \resizebox{55mm}{!}{\includegraphics{picture1.eps}}
\end{proof}

\begin{proof}[Proof of Theorem \ref{thm:1}]
We have seen using the lemma of homogeneity that it is sufficient
to prove Theorem \ref{thm:1} under the assumption $S(f) = S(g)$.
By Lemma \ref{lem:1} there exists a diffeomorphism $d$ such that
$(f \circ d)|F^2 \sim g|F^2$ and $S(f \circ d) = S(f)$. Now
applying Lemma \ref{lem:2} to $h = f \circ d$ and $g$ we get $f
\circ d \sim g$. But $f \circ d_t$ is a regular homotopy between
$f$ and $f \circ d$ proving that $f \sim g$.

It remains to show that the above regular homotopy  $H$ joining
$f$ and $g$ can be chosen in such a way that it defines a
prescribed bijection $i \colon S(f) \To S(g)$, that is $i_H = i$.
Recall that the bijection $i_H \colon S(f) \To S(g)$ depends only
on the choice of the diffeotopy $\varphi_t$ mentioned before Lemma
\ref{lem:1}, that clearly might induce any prescribed bijection
between $S(f)$ and $S(g)$. If $M^2$ is orientable, then the
diffeomorphism $d \colon M^2 \To M^2$ of Lemma \ref{lem:1} swaps
the singular points $p_1$ and $p_2$ (i.e. $d(p_1) = p_2$, $d(p_2)
= p_1$ and $d(p_i) = p_i$ for $i > 2$), and in the non-orientable
case $d(p_i) = p_i$ for $1 \le i \le k$. Finally the homotopy
constructed in Lemma \ref{lem:2} between $f \circ d$ and $g$ is a
singularity fixing homotopy in the sense of Definition
\ref{defn:2}. This completes the proof of Theorem \ref{thm:1}.
\end{proof}

The converse of Lemma \ref{lem:2} is true only in the following
form:

\begin{prop}
Suppose that $g$ and $h$ are locally generic mappings of $M^2$
into $\Real ^3$ such that $S(g)=S(h)$. Denote by $F^2$ the
complement of a small open neighborhood of the common singular
set. Then $g \sim h$ implies that $\imekv{g|F^2}{h|F^2}$.
\end{prop}

\begin{proof}
Let $H$ be a regular homotopy connecting $g$ and $h$. Then
\[S(H_t)= \{\,p_1(t), \dots, p_k(t)\,\},\] where $p_i(t)$ is a
smooth curve in $M^2$. The lemma of homogeneity gives a diffeotopy
$\varphi_t$ of $M^2$ such that $\varphi_0 = id_{M^2}$ and
$\varphi_t(S(H_0))=S(H_t)$ for every $t \in [0,1]$. The homotopy
$G_t = H_t \circ \varphi_t$ has the property that $S(G_t)=S(G_0)$
for $t \in [0,1]$ and connects $g$ with $h \circ \varphi_1$. Since
$\varphi_1$ permutes $S(g)$ (because $\varphi_1(S(g)) =
\varphi_1(S(H_0)) = S(H_1) = S(g)$) we can choose $\varphi_1$ to
map $F^2$ onto itself. $G_t|F^2$ is a regular homotopy between the
immersions $g|F^2$ and $(h|F^2) \circ (\varphi_1 | F^2)$ which
means by definition that $\imekv{g|F^2}{h|F^2}$.
\end{proof}

\begin{rem}
Modify Definition \ref{defn:1} of regular homotopy the following
way:
\begin{defn} \label{defn:2}
Locally generic mappings $f, g \colon M^2 \To \Real ^3$ are
\emph{regularly homotopic through a singularity fixing homotopy}
-- notation $f \sim_s g$ -- if $S(f) = S(g)$ and there exists a
smooth mapping $H \colon M^2 \times [0,1] \to \Real ^3$ such that
$H_0 = f$ and $H_1 = g$ and for every $t \in [0,1]$ the mapping
$H_t$ is locally generic with $S(H_t)=S(f)$. (That is the singular
points are kept fixed.)
\end{defn}

This gives a modification of the definition of image-homotopic
maps:

\begin{defn} Locally generic mappings $f, g \colon M^2 \To
\Real^3$ are \emph{image homotopic through a singularity fixing
homotopy} -- notation $\im(f) \sim_s \im(g)$ -- if $S(f) = S(g)$
and there is a $d: M^2 \To M^2$ diffeomorphism such that $f \circ
d \sim_s g$. Note that $d(S(f)) = S(g) = S(f)$ and $d$ can permute
the points of $S(f)$.
\end{defn}

Suppose that $S(f)$ or $S(g)$ is non-empty. The arguments above
show that $\im(f) \sim_s \im(g)$ if and only if $\abs{S(f)} =
\abs{S(g)}$. To prove this we only have to use diffeomorphisms
instead of diffeotopies since Lemma \ref{lem:2} remains true using
the new definition.

On the other hand $f \sim_s g$ implies that the immersions $f|(M^2
\setminus S(f))$ and $g|(M^2 \setminus S(f))$ are regularly
homotopic. But there are locally generic mappings $f, g \colon M^2
\To \Real ^3$ satisfying $S(f)=S(g)$ such that $f \sim g$ but
$f|(M^2 \setminus S(f)) \nsim g|(M^2 \setminus S(f))$. Take for
example $M^2 = RP^2$ and choose two arbitrary points $p, q \in
RP^2$. Denote $RP^2 \setminus \{p, q\}$ by $F^2$. Using the
notations of Lemma \ref{lem:1} we have that $H_1(F^2 ;
\mathbb{Z}_2) = \langle f_1, c_1 \rangle$. Define the cohomology
classes $\alpha, \beta \in H^1(F^2 ; \mathbb{Z}_2)$ by the
equalities $\alpha(f_1) = 0$, $\beta(f_1)=1$ and $\alpha(c_1) =
\beta(c_1) = 1$.  Then there exist locally generic mappings $f,g
\colon RP^2 \To \Real^3$ satisfying $S(f) = S(g) = \{p, q\}$ such
that $f|F^2$ and $g|F^2$ correspond to $\alpha$ and $\beta$ using
the bijection of Theorem \ref{thm:3} : Denote by $U$ the locally
generic mapping of $S^2$ to $\Real^3$ with singular points $p$ and
$q$ (this is unique up to singularity fixing homotopy). $B$ is the
famous Boy surface, $\overline{B}$ is the mirror image of $B$ (see
\cite{Pinkall}). Then the connected sums $f = B \# U$ and $g =
\overline{B} \# U$ satisfy the above conditions. Clearly $f|F^2
\nsim g|F^2$, thus $f \nsim_s g$. This provides examples of
locally generic mappings $f$ and $g$ such that $\im(f) \sim_s
\im(g)$, but $f \nsim_s g$.
\end{rem}

\section{Projections of regular homotopies}

Suppose that $M^2$ is a closed  connected surface. In the previous
sections we examined the path-components of the space of locally
generic mappings of $M^2$ into $\Real^3$ endowed with the
$C^{\infty}$ topology. This space, denoted by $\lgen(M^2,
\Real^3)$, is closely connected with the space $\imm(M^2,
\Real^4)$ of immersions of $M^2$ into $\Real^4$ (also considered
with the $C^{\infty}$ topology). This connection is realized by
projections of $\Real^4$ onto $\Real^3$.

\begin{defn}
A mapping $\pi \colon \Real^k \To \Real^l$ $(k > l)$ is a
projection if it is linear and surjective. Denote the space of all
projections from $\Real^k$ to $\Real^l$ by $\proj(\Real^k,
\Real^l)$.
\end{defn}

Our starting point is the following result of Mather
\cite{Mather}:

\begin{prop} \label{prop:3}
Suppose that $F \in \imm(M^2, \Real^4)$ is an immersion. Then for
almost every $\pi \in \proj(\Real^4, \Real^3)$ (in the sense of
Lebesgue measure) the mapping $\pi \circ F$ is locally generic.
\end{prop}

For every projection $\pi \in \proj(\Real^4, \Real^3)$ let
\[\imm_{\pi}(M^2, \Real^4) = \left\{\, F \in \imm(M^2, \Real^4) \colon \pi \circ F \in
\lgen(M^2, \Real^3)\,\right\}\] be the subspace of $\imm(M^2,
\Real^4)$. There is a natural mapping \[p_{\pi} \colon
\imm_{\pi}(M^2, \Real^4) \To \lgen(M^2, \Real^3)\] defined by the
formula $p_{\pi}(F) = \pi \circ F$ for every $F \in
\imm_{\pi}(M^2, \Real^4)$.

In this section our aim is to examine the path-components of
$\imm_{\pi}(M^2, \Real^4)$.

\begin{defn}
Two immersions $F,G \in \imm_{\pi}(M^2, \mathbb{R}^4)$ are called
\emph{$\pi$-homotopic} (denoted by $F \sim_{\pi} G$) if they are
in the same path-component of $\imm_{\pi}(M^2, \mathbb{R}^4)$.
\end{defn}

First let us recall that two immersions $F, G \in \imm(M^2,
\Real^4)$ are regularly homotopic if and only if $e(F) = e(G)$,
where $e(F)$ denotes the (twisted) Euler-number of the normal
bundle of the immersion $F$. It is clear that if $F \sim_{\pi} G$
then $e(F) = e(G)$ and $\abs{S(\pi \circ F)} = \abs{S(\pi \circ
G)}$. We are going to prove that if $S(\pi \circ F)$ (and $S(\pi
\circ G)$) are non-empty then the converse also holds.

Suppose that $F, G \colon M^2 \To \Real^4$ are immersions. From
Proposition \ref{prop:3} it is clear that for almost every
projection $\pi \in \proj(\Real^4, \Real^3)$ both $f = \pi \circ
F$ and $g = \pi \circ G$ are locally generic. Fix such a
projection $\pi$.

We are going to define the sign of every point in $S(f)$ (and in
$S(g)$). (An equivalent definition can be found in \cite{Saeki}).
Take a cross-cap point $p \in S(f)$. Choose orientations of
$\Real^4$ and of $\Real^3$. We will define the sign of $p$ as
follows.

Fix local coordinates $(x_1, x_2)$ on a neighborhood $U_p$ of $p$
and $(y_1, y_2, y_3)$ centered at $f(p)$ such that $f$ has the
following normal form: \[y_1 \circ f = x_1^2, \,\,\, y_2 \circ f =
x_2,\,\,\, y_3 \circ f = x_1x_2.\] Suppose that $U_p$ is so small
that $F|U_p$ is an embedding. The sign of $p$ will depend only on
$F|U_p$ (thus the definition is local). Set $D_{\eps} = \{y_1^2 +
y_2^2 + y_3^2 \le \eps \} \subset \Real^{3}$ and choose $\eps > 0$
sufficiently small such that $\widetilde{D}_{\eps} =
f^{-1}(D_{\eps}) \subset U_p$. The set $\widetilde{D}_{\eps}$ is a
closed disc neighborhood of $p$ in $M^2$ and $\partial
\widetilde{D}_{\eps} = f^{-1}(\partial D_{\eps})$ (see Lemma 2.2
in \cite{Saeki}). Let $L$ be the closure of the double point set
of $f|\widetilde{D}_{\eps}, i.e. $ $L = \{\, x \in
\widetilde{D}_{\eps} \colon f^{-1}(f(x)) \cap \widetilde{D}_{\eps}
\neq \{x\} \,\} \cup \{p\}$. Then $L$ is a one-dimensional smooth
submanifold of $\widetilde{D}_{\eps}$ and $L \cap
\partial\widetilde{D}_{\eps}$ consists of two points $p_1$
and $p_2$. We fix an orientation of $\partial
\widetilde{D}_{\eps}$ and take an oriented base $(u)$ (resp.
$(v)$) of the tangent space $T_{p_1}(\partial
\widetilde{D}_{\eps})$ (resp. $T_{p_2}(\partial
\widetilde{D}_{\eps})$). Then $(df_{p_1}(u), df_{p_2}(v))$ is a
base of $T_q(\partial D_{\eps})$, where $q = f(p_1) = f(p_2)$. We
may assume that $(df_{p_1}(u), df_{p_2}(v), \xi)$ is a positive
basis of $T_q \Real^3$, where $\xi$ is the outward normal vector
of $\partial D_{\eps}$, exchanging $p_1$ and $p_2$ if necessary.
Now orient $L$ from $p_2$ to $p_1$.

Denote by $\nu$ a positive basis vector of $T_pL$. (If $M^2$
possess a Riemannian metric, then choose $\nu$ to be a
unit-vector. This way $\nu$ is unique up to the orientation of
$\Real^3$.) Orient $\ker \pi$ in such a way that together with the
orientation of $\im \pi =\Real^3$ we obtain the fixed orientation
of $\Real^4 = \im \pi \oplus \ker \pi$. Using this direct sum
decomposition of $\Real^4$ the mapping $F \colon M^2 \To \Real^4$
can be written in the form $F = (\pi \circ F, F^*) = (f, F^*)$.
Since $F$ is an immersion at the point $p$ and $p \in S(f)$ , we
have $dF^*(\nu) \neq 0$. After all this preparation we can now
define the sign of $p$.

\begin{defn} \label{defn:4}
The cross-cap point $p$ is \emph{positive} if $dF^*(\nu) > 0$ and
is \emph{negative} if $dF^*(\nu) <0$. We denote by $p(F)$ (resp.
$n(F)$) the number of positive (resp. negative) cross-cap
singularities of the locally generic mapping $f = \pi \circ F
\colon M^2 \To \Real^3$.
\end{defn}

The following proposition is  a special case of Proposition 2.5 in
\cite{Saeki}.

\begin{prop}
Suppose that $\Real^4$ and $\Real^3$ are oriented. Then we always
have
\[e(F) = p(F) - n(F),\] where $e(F) \in \mathbb{Z}$ is the
(twisted) Euler-number of the normal bundle of the immersion $F$.
\end{prop}

Thus the immersions $F$ and $G$ are regularly homotopic if and
only if \[p(F) - n(F) = e(F) = e(G) = p(G) - n(G).\] On the other
hand Theorem \ref{thm:1} states that if $S(f)$ and $S(g)$ are
non-empty then $f \sim g$ if and only if \[p(F) + n(F) =
\abs{S(f)} = \abs{S(g)} = p(G) + n(G).\] Comparing the preceding
two chains of equations we have that if both $f$ and $g$ are
singular then $[F \sim G$ and $f \sim g] \Leftrightarrow [p(F) =
p(G)$ and $n(F) = n(G)]$. The following theorem implies that in
this case we can even find a regular homotopy between $F$ and $G$
whose projection is a regular homotopy between $f$ and $g$, i.e.
$F \sim_{\pi} G$.

\begin{thm} \label{thm:4}
Suppose that $M^2$ is a closed connected surface, $F, G \colon M^2
\To \Real^4$ are immersions and $\pi \colon \Real^4 \To \Real^3$
is a projection such that $f = \pi \circ F$ and $g = \pi \circ G$
are both locally generic and singular. Then the following are
equivalent:

$(1)$ There exists a regular homotopy $H \colon M^2 \times [0,1]
\To \Real^4$ between $F$ and $G$ such that $\widetilde{H} = \pi
\circ H$ is a regular homotopy between $f$ and $g$, i.e. $F
\sim_{\pi} G$.

$(2)$ The numbers of positive and negative cross-caps of $f$ and
$g$ are the same, i.e. $p(F) = p(G)$ and $n(F) = n(G)$.
\end{thm}

\begin{proof}
First we prove the implication (1) $\Rightarrow$ (2). In this case
$f \sim g$, thus using Proposition \ref{prop:1} we have that
$\abs{S(f)} = \abs{S(g)}$. From Definition \ref{defn:4} it is
clear that the signs of the singular points do not change during a
regular homotopy. (Here we did not use the assumption that
$\abs{S(f)}>0, \abs{S(g)} > 0$.)

Now we are going to prove the implication (2) $\Rightarrow$ (1).
Since $p(F) = p(G)$ and $n(F) = n(G)$ there exists a bijection $i
\colon S(f) \To S(g)$ that preserves the signs of the cross-cap
points. By Theorem \ref{thm:1} there is a regular homotopy
$\widetilde{H}$ between $f$ and $g$ such that $i_{\widetilde{H}} =
i$. We shall construct a regular homotopy $H$ between the
immersions $F$ and $G$ such that $\pi \circ H = \widetilde{H}$ as
follows: Choose an arbitrary Riemannian metric on $M^2$. In the
paragraph preceding Definition \ref{defn:4} we saw that in this
case for every $t \in [0,1]$ and $p \in S(\widetilde{H}_t)$ there
is a unique positive unit-vector $\nu_t(p) \in T_pM^2$ tangent to
the double-point curve $L$ crossing $p$. The singular sets
$S(\widetilde{H}_t)$ for $t \in [0,1]$ define curves $\gamma_1,
\dots, \gamma_k$ on $M^2$ such that $S(H_t) = \left\{ \,
\gamma_j(t) \colon 1 \le j \le k \, \right\}$ for every $t \in
[0,1]$. The points $\gamma_i(0)$ and $\gamma_i(1 )$ have the same
sign ($1 \le i \le k$). Suppose for example that for a fix $i$
both $\gamma_i(0)$ and $\gamma_i(1)$ are positive cross-cap
points. Introduce the notation $\nu_i(t) = \nu_t(\gamma_i(t))$,
then $dF^*(\nu_i(0)) > 0$ and $dG^*(\nu_i(1)) > 0$. (Here $F^*$
denotes the fourth coordinate function of $F$ in $\mathbb{R}^4$.)
Using the Levi-Civita connection of the Riemannian manifold $M^2$
we may consider the exponential mapping on $M^2$. Since $[0,1]$ is
compact, there exists $\eps > 0$ such that for every $t \in [0,
1]$ the mapping $h \colon [-\eps,\eps] \times [0,1] \To M^2$
satisfying \[h(s,t) = \exp_{\gamma_i(t)}(s \cdot \nu_i(t))\] is
defined and for every $t \in [0,1]$ the mapping $h_t(s) = h(s,t)$
is an embedding of $[-\eps,\eps]$ into $M^2$ ($h_t$ is a geodetic
curve). Define the function $H_t^*$ on $\im h_t$ using the
following formula:
\[H_t^*(h_t(s)) = (1 - t) \cdot F^*(h_0(s)) +
t \cdot G^*(h_1(s))\] for $s \in [-\eps, \eps]$. Note that $H_0^*|
\im h_0 = F^*| \im h_0$ and $H_1^*| \im h_1 = G^*| \im h_1$. Thus
we can extend $H^*$ to an open neighborhood of $\bigcup_{t \in
[0,1]}\left(\im h_t \times \{t\} \right)$ in $M^2 \times [0,1]$ as
a smooth function. From the construction of $H_t^*$ it is clear
that
\[dH_t^*(\nu_i(t)) = (1-t) \cdot dF^*(\nu_i(0)) + t \cdot
dG^*(\nu_i(1)) > 0,\] which implies that the mapping $H_t =
(\widetilde{H}_t, H_t^*)$ is an immersion at the point
$\gamma_i(t)$ for every $t \in [0,1]$. Repeat the preceding
extension process for every $1 \le i \le k$ and afterwards extend
the obtained $H^*$ to the whole cylinder $M^2 \times [0,1]$ in
such a way that $H_0^* = F^*$ and $H_1^* = G^*$. The mapping $H =
(\widetilde{H}, H^*)$ is a regular homotopy connecting $F$ and $G$
whose projection is $\widetilde{H}$.
\end{proof}

Putting together our previous results we obtain the following
theorem:

\begin{thm}
If $F, G \in \imm_{\pi}(M^2, \mathbb{R}^4)$ then \[ F \sim_{\pi} G
\Leftrightarrow [F \sim G \,\, \text{and} \,\, \pi \circ F \sim
\pi \circ G].
\]
\end{thm}

\begin{proof}
First we suppose that $S(\pi \circ F)$ or $S(\pi \circ G)$ is
non-empty. Theorem \ref{thm:4} states that $F \sim_{\pi} G
\Leftrightarrow [p(F) = p(G)$ and $n(F) = n(G)]$. We have seen in
the paragraph preceding Theorem \ref{thm:4} that $ [p(F) = p(G)$
and $n(F) = n(G)] \Leftrightarrow [F \sim G$ and $\pi \circ F \sim
\pi \circ G]$.

Now we consider the case when both $S(\pi \circ F)$ and $S(\pi
\circ G)$ are empty, i.e. $f = \pi \circ F$ and $g = \pi \circ G$
are immersions. If $f \sim g$ then any regular homotopy connecting
$f$ and $g$ can be lifted to a regular homotopy between $F$ and
$G$, thus $F \sim_{\pi} G$. This proves the implication $F
\sim_{\pi} G \Leftarrow [F \sim G$ and $\pi \circ F \sim \pi \circ
G]$. The other implication is trivial.
\end{proof}
% ----------------------------------------------------------------
\bibliographystyle{amsplain}
\bibliography{topology}
\end{document}